\documentclass[12pt]{article}
\usepackage{amsmath,amssymb,euscript}
\newtheorem{theorem}{Theorem}

\newtheorem{proposition}[theorem]{Proposition}
\newtheorem{corollary}[theorem]{Corollary}
\newtheorem{definition}[theorem]{Definition}

\begin{document}

\title{\bf The Rokhlin lemma for homeomorphisms of a Cantor set}
\author{{\bf S.~Bezuglyi}\\
Institute for Low Temperature Physics, Kharkov, Ukraine\\ \and{\bf
A.H.~Dooley}\\ University of New South Wales, Sydney, Australia\\ \and {\bf
K.~Medynets}\\ Institute for Low Temperature Physics, Kharkov, Ukraine}

\maketitle

\newcommand{\e}{\varepsilon}
\newcommand{\N}{{\mathbb N}}
\newcommand{\h}{Homeo(X)}
\newcommand{\M}{\mathcal{M}_1(X)}
\newcommand{\la}{\lambda}
\newcommand{\per}{{\cal P}er}
\newcommand{\ap}{{\cal A}p}

\begin{abstract}
For a Cantor set $X$, let $Homeo(X)$ denote the group of all homeomorphisms
of $X$. The main result of this note is the following theorem. Let $T\in
Homeo(X)$ be an aperiodic homeomorphism, let $\mu_1,\mu_2,\ldots,\mu_k$ be
Borel probability measures on $X$, $\e> 0$, and $n\ge 2$. Then there exists
a clopen set $E\subset X$ such that the sets $E,TE,\ldots, T^{n-1}E$ are
disjoint and $\mu_i(E\cup TE\cup\cdots\cup T^{n-1}E) > 1 - \e,\ i=
1,\ldots,k$. Several corollaries of this result are given. In particular, it
is proved that for any aperiodic $T\in \h$ the set of all homeomorphisms
conjugate to $T$ is dense in the set of aperiodic homeomorphisms.
\end{abstract}

\noindent {\bf 0. Introduction}. One of the most useful results in ergodic
theory which has many important applications is the Rokhlin lemma [R]. This
statement asserts that given an aperiodic (non-singular) automorphism $T$ of
a standard measure space $(X,\mathcal{B},\mu)$, $\e>0$, and $n\ge 2$, there
exists a measurable subset $E\subset X$ such that $E,TE,\ldots,T^{n-1}E$ are
pairwise disjoint and $\mu(E\cup TE\cup\ldots \cup T^{n-1}E)>1-\e$. It
immediately follows from this result that the set of periodic automorphisms
is dense in the group of all nonsingular automorphisms of
$(X,\mathcal{B},\mu)$ with respect to the metric $d(T_1,T_2) = \mu\{(x\in X
: T_1x \neq T_2x)\}$ where $T_1,T_2\in Aut(X,\mathcal{B},\mu)$.
Subsequently, the Rokhlin lemma was generalized in various directions (see,
for example, [AP, EP, FL, LW, OW]). It is well known that the Rokhlin lemma
is also related to amenability.

Our goal is to prove a version of the Rokhlin lemma in the context
of Cantor dynamics. To start with, we need to consider a topology
on $\h$ analogous to the metric $d$. It is well known that $\h$
equipped with the topology $\tau_w$ of uniform convergence (the
$\sup$-topology) is a Polish space. But $\tau_w$ is too weak to
trace the dynamical properties of homeomorphisms. For instance,
the set of pointwise periodic homeomorphisms is not dense in $\h$
with respect to $\tau_w$ (more details see in [BDK1, BDK2]). In
[BK1, BK2, BDK1, BDK2], we introduced and studied a topology
$\tau$ which has its origin in measurable dynamics. It was shown
in [BDK1, BDK2] that $\tau$ plays the same role in Borel and
Cantor dynamics as the metric $d$ plays in ergodic theory. Recall
the definition of the uniform topology $\tau$. For $T,S\in \h$,
define $E(S,T)=\{x\in X : Tx\neq Sx\}\cup\{x\in X : S^{-1}x\neq
T^{-1}x\}$. Let $\mathcal{M}_1(X)$ denote the set of all Borel
probability measures on $X$.

\begin{definition} The uniform topology $\tau$ is defined by the
base of neighborhoods $\mathcal{U}=\{U(T;\mu_1,\mu_2,\ldots,\mu_n)\}$ where
$$
U(T,\mu_1,\mu_2,\ldots,\mu_n,\e)=\{S\in Homeo(X)\,|\,
\mu_i(E(S,T))<\e,\;i=1,\ldots,n\},
$$
$T\in Homeo(X)$, $\mu_1,\mu_2,\ldots,\mu_n\in\mathcal{M}_1(X)$, and $\e>0$.
\end{definition}

$\h$ is a Hausdorff topological group with respect to $\tau$.

Remark that if one considers the set $E_0(T,S) = \{x \in X : Sx \neq Tx\}$
instead of $E(T,S)$ in the definition of $\tau$, then the topology defined
by $E_0(T,S)$ is equivalent to $\tau$.

The main result of this paper is the Rokhlin lemma proved for an arbitrary
aperiodic homeomorphism of Cantor dynamics. Let $\ap$ denote the set of all
aperiodic homeomorphisms of $X$.

\begin{theorem}\label{Rokhlin}
Let $T\in \ap$, $\mu_1,\ldots,\mu_k \in \M$, $\e>0$, and let $n\ge 2$. Then
there exists a clopen set $E\subset X$ such that sets $E,TE,\ldots,T^{n-1}E$
are pairwise disjoint and
$$
\mu_i(E\cup TE\cup\ldots \cup T^{n-1}E)>1-\e\ \ \mbox{ for all }1\leq i\leq
k.
$$
\end{theorem}

In contrast to measurable dynamics, we do not assume that $\mu_i$ is
invariant (or non-singular) with respect to $T$. Furthermore, we prove in
Theorem \ref{Rokhlin} that the base $E$ of the tower $(E, TE,...,T^{n-1}E)$
can be chosen clopen.

The proof of Theorem \ref{Rokhlin} is based on the following result. It
turns out that a Cantor set $X$ can be represented as a disjoint union of
clopen $T$-towers of sufficiently large height where $T$ is an aperiodic
homeomorphism. Then we can reconstruct these towers to get new ones such
that the union of bases and tops is of small measure.

We prove several statements which are immediate consequences of Theorem
\ref{Rokhlin}. In particular, we will show that every aperiodic
homeomorphisms may be approximated in $\tau$ by periodic homeomorphisms with
finite periods. We will further prove that, given $S\in \ap$, the set
$\{TST^{-1} : T\in \h\}$ is $\tau$-dense in $\ap$. This property, called the
Rokhlin property, was also studied in [GK, GW, R].
\\

\noindent {\bf 1. $T$-towers covering a Cantor set}. Let $T\in \h$ and let
$E$ be a clopen set such that $E\cap T^i(E) =\emptyset, i=1,...,n-1$. Then
$\xi = (E,TE,...,T^{n-1}E)$ is called a $T$-tower with base $E$ and top
$T^{n-1}(E)$. The number $h(\xi)= n$ is the height of $\xi$ and the clopen
set $C = \bigcup_{i=0}^{n-1}T^iE$ is called the support of $\xi$.

If $T$ is a minimal homeomorphism of a Cantor set $X$, then for every clopen
set $E$ one can easily find a finite collection of disjoint $T$-towers $\Xi
= (\xi_1,...,\xi_m)$ such that their supports $C_i,\ i=1,...,m$, cover $X$
and $E$ is the union of bases of these towers. Note that $T$ maps the union
of tops of $\xi_i$ onto $E$.

In this section we show that a similar cover exists for any aperiodic
homeomorphism $T$ of $X$. The following proposition will be used in our
proof of the Rokhlin lemma.

\begin{proposition}\label{cover}
Let $T$ be an aperiodic homeomorphism of a Cantor set $X$. Given a
positive integer $n\ge 2$, there exists a partition of $X$ into a
finite number of $T$-towers $\Xi=\{\xi_1,\ldots,\xi_m\}$ such that
the height $h(\xi_i)$ of every tower is at least $n$.
\end{proposition}
{\it Proof.} Since $T$ is aperiodic, every point $x\in X$ has a clopen
neighborhood $V_x$ such that the sets $V_x,TV_x,\ldots,T^{n-1}V_x$ are
pairwise disjoint. The open cover $\{V_x\,:\,x\in X\}$ of $X$ contains a
finite subcover $\{V_1,V_2,\ldots,V_k\}$. Setting
$$\begin{array}{lll}
U_1=V_1\\ U_2=V_2-U_1\\ \dots\dots\dots\dots
\\U_k=V_k-(U_1\cup U_2\cup\ldots\cup U_{k-1})
\end{array}
$$
we get a disjoint finite cover of $X$ by clopen sets such that $U_i\cap
T^jU_i = \emptyset,\ j=1,...,n-1,\ i=1,...,k$. Throughout the proof we will
deal with clopen sets only.
\medskip

(a) The set $\Xi$ will be constructed inductively. At the first step, we
start with the $T$-tower $\xi_1 = (U_1, TU_1,...,T^{n-1}U_1)$.

Let $C_1=\bigcup_{s=0}^{n-1}T^{s}U_1$ be the support of $\xi_1$. If $C_1 =
X$, we are done. If $C_1$ is a proper subset of $X$, consider the sets
$U_i^1=U_i-C_1,\ i=2,3,\ldots,k$. Note that some of the sets $U_i^1$ may be
empty. Without loss of generality assume that $U_2^1 \neq \emptyset$. Notice
that the first set from $C_1$ which could meet the $T$-orbit of $U_2^1$ is
$U_1$. Define
$$
U_2^1(i)=\{x\in U_2^1\,:\, T^ix\in U_1,\ T^jx\notin U_1,\ 0\leq j\leq
i-1\},\ i=1,\ldots,n-1,
$$
and
$$
U_2^1(0)=\{x\in U_2^1\,:\,T^jx\notin U_1,\mbox{ for all }1\leq j\leq n-1\}.
$$
In particular, some of the sets $U_2^1(i),\ i=1,...,n-1$, may be empty. Each
set $U_2^1(i)$ is the base of the $T$-tower
$$
\xi_2^1(i)=\{U_2^1(i),TU_2^1(i),\ldots,T^{n-1+i}U_2^1(i)\},
\mbox{ for all }i=0,1,\ldots,n-1.
$$
We see that the first $i-1$ sets of the tower $\xi_2^1(i)$ do not intersect
$\xi_1$ and the other sets form a subtower in $\xi_1$. Remark that the
towers $\xi_2^1(i)$ and $\xi_2^1(j)$ are obviously disjoint whenever $i\neq
j$ but each tower $\xi_2^1(i),\ 1\leq i\leq n-1$, meets $\xi_1.$ To make all
of them disjoint, we delete all sets $T^iU_2^1(i)$, $1\leq i\leq n-1$, from
the base of the tower $\xi_1$. Take the set $U_1^1 = U_1 -
\bigcup\limits_{i=1}^{n-1}T^iU_2^1(i)$ as the base of a subtower of $\xi_1$.
Denote this subtower by $\xi_1^1$. At the end of the first step, the
collection of disjoint $T$-towers $\Xi(1) =
\{\xi_1^1,\xi_2^1(0),\xi_2^1(1),\ldots,\xi_2^1(n-1)\}$ each of height at
least $n$. Denote by $C_1^1, C_2^1(0), C_2^1(1),...,C_2^1(n-1)$ the supports
of corresponding towers from $\Xi(1)$.
\medskip

(b) We note that $C_1^1\cup C_2^1(0)\cup C_2^1(1)\cup\cdots \cup
C_2^1(n-1) = C_1 \cup C_2^1$ where
$C_2^1=\bigcup\limits_{s=0}^{n-1}T^i U_2^1$. If $C_1\cup C_2^1$
covers $X$, the proof is complete. If not, define the sets
$U_i^2=U_i^1-C_2^1,$ for all $i=3,4,\ldots,k$. It is easily seen
that $\{C_1^1, C_2^1(0),
C_2^1(1),...,C_2^1(n-1),U_3^2,...,U_k^2\}$ is a partition of $X$.
Without loss of generality we assume that $U_3^2\neq \emptyset$.
For each tower $\xi$ of $\Xi(1)$, put
$$U_3^2(\xi)=\{x\in U_3^2 : T^ix\mbox{ meets }\xi\mbox{ for some }1\leq i\leq n-1\}.$$
And let $U_3^2(0)=U_3^2-\bigcup_{\xi\in\Xi(1)}U_3^2(\xi)$. Put
$\eta=\{U_3^2(0),\ldots, T^{n-1}U_3^2(0)\}$. Notice that $\eta$
does not intersect any $\xi\in\Xi(1)$. We can apply now the
construction used in case (a) to the set $U_3^2(\xi)$ and the
tower $\xi\in\Xi(1)$. This means that we can partition
$U_3^2(\xi)$ into the sets $\{U_3^2(\xi,i)\,:\,1\leq i\leq n-1\},$
where $U_3^2(\xi,i)$ contains all points whose orbits meet $\xi$
precisely at the $i$-th step. Thus, we have defined a $T$-tower
$\xi(i),\ i=1,...,n-1,$ with base $U_3^2(\xi,i)$. The height of
this tower is $h(\xi)+i$ where $h(\xi)$ is the height of $\xi$.
The towers $\xi(i)$ remove at most $n-1$ subtowers from $\xi$.
Denote by $\xi'$ the remaining part of $\xi$. Thus, $\xi$ can be
represented as a disjoint union of at most $n$ $T$-towers. At the
end of this step, the set $\Xi(2)$ consists of all $\xi', \xi(i)$
and $\eta$ where $1\le i < n$ and $\xi$ runs over $\Xi(1)$.
Clearly, we have constructed a disjoint set of $T$-towers of
heights at least $n$. The union of supports of these towers equals
$C_1\cup C_2^1\cup C_3^2$ where
$C_3^2=\bigcup\limits_{i=0}^{n-1}T^iU_3^2$.
\medskip

(c) If $C_1\cup C_2^1\cup C_3^2 = X$, the proof is finished. If $C_1\cup
C_2^1\cup C_3^2$ is a proper subset of $X$ we consider the sets
$U_j^3=U_j^2-C_3^2$, for all $j=4,5,\ldots,k,$ and apply the above
construction again. It is obvious that that after at most $k-1$ steps we
will have constructed a collection of towers $\Xi$ satisfying the conditions
of the proposition. \hfill$\square$

\begin{corollary}\label{Section} Let $T$ be an aperiodic homeomorphism of
$X$. Suppose that there exist a clopen set $F$ and $N\in \N$ such that
$X=\bigcup\limits_{i=0}^N T^iF$. Then for any integer $n \geq N$ there
exists a finite partition $\Xi=\{\xi_0,\ldots,\xi_m\}$ of $X$ into
$T$-towers of height at least $n$ such that the union of the bases of these
towers is a subset of $F$.
\end{corollary}
{\it Proof.} This result can be proved by the method used in Proposition
\ref{cover}. To do this, we start with a finite cover $\{U_1,...,U_k\}$ of
$F$ by disjoint clopen sets such that $U_i\cap T^jU_i = \emptyset,\ 0< j
<n,\ 1\le i\le k$. Then we can construct disjoint clopen $T$-towers of
height at least $n$ with bases containing in the $U_i$'s as above. Since
$X=\bigcup_{i=0}^NT^iF$ and $n\geq N$ we that the towers we have constructed
cover $X$. \hfill$\square$
\\

\noindent {\bf 2. Proof of the Rokhlin lemma}. This proof will be based on
the decomposition of $X$ into a finite set of towers found in Proposition
\ref{cover}. In this proof, we again work with clopen sets only.
\\

\noindent {\it Proof of Theorem $2$}. Choose $m\in \N$ such that
$1/m<\e$. Using Proposition \ref{cover} find a finite set of
disjoint $T$-towers $\Lambda =(\la_1,...,\la_q)$, consisting of
clopen sets, such that $h(\la_i)\ge nm^k,\ i=1,...,q$. Denote by
$V_1,\ldots,V_q$ and $B_1,...,B_q$ the tops and bases of towers
$\la_1,\ldots,\la_q$, respectively. It follows from the
construction of towers that
$$
T(\bigcup_{i=1}^q V_i) = \bigcup_{i=1}^q B_i.
$$

Fix the measure $\mu_1$. Let
$$
F_l=\bigcup\limits_{s=1}^q\ \bigcup\limits_{j=(l-1)nm^{k-1}}^{lnm^{k-1}-1}
T^{-j}V_s,\ \ l=1,\ldots,m.
$$
Then $F_1,F_2,\ldots,F_m$ are disjoint clopen sets. Amongst these sets there
exists at least one of them, say $F_{l_1}$, whose $\mu_1$-measure is not
greater than $1/m$. Denote this set by $F^{l_1},\ 1\leq l_1\leq q$. By
construction, $F^{l_1}$ consists of $q$ towers
$$
\la_s^{l_1}=\{T^{-j}V_s\;:\;(l_1-1)nm^{k-1}\leq j\leq l_1nm^{k-1}-1 \},\ \
s=1,2,\ldots,q,
$$
such that $h(\la_s^{l_1}) \ge nm^{k-1}$. Let $V_s^{l_1}$ denote the top of
the tower $\la_s^{l_1}$, $s=1,2,\ldots,q$.

Define the sets
$$
F_l^{l_1}=\bigcup\limits_{s=1}^q
\bigcup\limits_{j=(l-1)nm^{k-2}}^{lnm^{k-2}-1}T^{-j}V_s^{l_1},\ \
l=1,2,\ldots,m.
$$
Since $F_1^{l_1},F_2^{l_1},\ldots,F_m^{l_1}$ are disjoint sets, at
least one of them has the $\mu_2$-measure not greater than $1/m$.
Denote this set by $F^{l_1,l_2}$. Since $F^{l_1,l_2}\subset
F^{l_1}$, we have that $\mu_i(F^{l_1,l_2}) \le 1/m,\ i=1,2$.

We continue this construction by induction. At the last step, we have found
a set $F^{l_1,\ldots,l_{k-1}}$ which consists of $q$-towers of height at
least $nm$. Denote by
$V_1^{l_1,\ldots,l_{k-1}},...,V_q^{l_1,\ldots,l_{k-1}}$ the tops of these
towers. Define
$$
F_l^{l_1,\ldots,l_{k-1}}= \bigcup\limits_{s=1}^q
\bigcup\limits_{j=(l-1)n}^{ln-1}T^{-j}V_s^{l_1,\ldots,l_{k-1}},\ \
l=1,2,\ldots,m.
$$
As above, at least one of the sets
$F_1^{l_1,\ldots,l_{k-1}},...,F_q^{l_1,\ldots,l_{k-1}}$ has $\mu_k$-measure
not greater than $1/m$. We denote this set by $F^{l_1,...,l_k} =F_0$. Thus,
by construction of $F_0$, we have that
$$
\mu_i(F_0)\le 1/m<\e,\ \ i=1,2,\ldots,k,
$$
and, moreover, $\bigcup\limits_{j=0}^\infty T^jF_0=X$ (actually, this union
is finite). We see that $F_0$ is formed by $q$ $T$-towers which have the
height $n$. Let $V_1^{l_1,...,l_k},...,V_q^{l_1,...,l_k}$ be the tops of
these towers. Denote
$$
F = \bigcup_{s=1}^q V_s^{l_1,...,l_k}.
$$
Then the above inequality for $\mu_i(F_0)$ implies that
$\mu_i(\bigcup\limits_{j=0}^{n-1}T^{-j}F)<\e,$ for all $i=1,2,\ldots,k.$
Note that there exists $K\in \N$ such that $T^K(V_s^{l_1,...,l_k}) = V_s$
for all $s = 1,...,q$.

Now we define a new finite collection $\Xi$ of $T$-towers whose bases are
clopen subsets of $F$. Let $D_{sp} = T^{-1}B_p\cap V_s$ and $C_{sp} =
T^{-K}D_{sp},\ 1\le s,p\le q$ (some of these sets may be empty). Clearly,
$\bigcup_{s,p=1}^{q}C_{sp} = F$. Then
$(C_{sp},TC_{sp},...,T^{h(sp)-1}C_{sp})$ constitutes a $T$-tower $\xi(sp)\in
\Xi$ where the height $h(sp)$ of $\xi(sp)$ is defined from the condition
$T^{h(sp)-1}C_{sp} \subset T^{-1}F,\ 1\leq s,p \leq q$. We get that the set
$\Xi= (\xi(sp) : 1\le s,p \le q)$ consists of $T$-towers whose supports form
a partition of $X$ and
$$
\mu_i(\bigcup\limits_{l=0}^{n-1}\bigcup\limits_{s,p=1}^q
T^{h(sp)-1-l}C_{sp})<\e, \ \ \ i=1,2,\ldots,k.
$$

To finish the proof, define the clopen set $E$:
$$
E=\bigcup\limits_{s,p=1}^q\bigcup\limits_{j= 0}^{L(sp)} T^{jn}C_{sp}
$$
where $L(sp) = [n^{-1}(h(sp)-1)]-1$. It is easily seen that $E\cap T^jE =
\emptyset,\ j=1,...,n-1$, and
$$
\mu_i(E\cup TE\cup\ldots \cup T^{n-1}E)>1-\e, \ \ 1\leq i\leq k.
$$
\hfill$\square$
\\

\noindent {\bf 3. Corollaries}. We give several immediate consequences of
the result proved above. The next two corollaries can be easily deduced from
the proof of Theorem \ref{Rokhlin}.

\begin{corollary}\label{C1} Suppose $T$ is an aperiodic homeomorphism of
$X$. Let $\mu_1,\ldots,\mu_k$ be Borel probability measures, $n\ge 2$ a
positive integer, and $\e>0.$ Then there exists a clopen partition of $X$
into a finite number of $T$-towers $\Xi=\{\xi_1,\ldots,\xi_q\}$ such that:
\begin{enumerate}
\item[(i)] for every tower $\xi\in\Xi$, $h(\xi) \geq n$; \item[(ii)] for $i=
1,...,k$,
$$
\mu_i(\bigcup\limits_{s=1}^q\bigcup\limits_{j=0}^{h(\xi_s)- n}T^jB_s)>1 -\e,
$$
\end{enumerate}
where $B_s$ is the base of $\xi_s$.
\end{corollary}

\begin{corollary} Suppose $T$ is an aperiodic homeomorphism of $X$. Let
$\mu_1,\ldots,\mu_k$ be Borel probability measures and let $\e>0$.
Then there exists a clopen set $A$ such that $\bigcup_{j\ge 0}T^jA
= X$ and $\mu_i(A) < \e,\ i=1,...,k$.
\end{corollary}

Let $T$ be an aperiodic homeomorphism of a Cantor set $X$. We recall the
definition of {\it the full group} $[T]$ of homeomorphisms generated by $T$.
By definition, a homeomorphism $\gamma\in Homeo(X)$ belongs to $[T]$ if the
sets $X_n = \{x\in X : \gamma x=T^nx\}$ form a partition of $X$. In other
words, there exists an integer-valued function $n(x)$ such that $\gamma
x=T^{n(x)}x$ for all $x\in X$.

Let $\per$ denote the set of pointwise periodic homeomorphisms of $\h$.
Consider the subset $\per_0\subset \per$ consisting of all homeomorphisms
with finite period, that is $P\in \per_0$ if and only if there exists $m\in
\N$ such that $P^mx =x$ for all $x\in X$. This means that $X$ can be
decomposed into a finite union of clopen sets $X_p$ such that the period of
$P$ at each point from $X_p$ is exactly $p$.

The following statement shows that every aperiodic homeomorphism can be
approximated by a periodic homeomorphism in the uniform topology. This
result has analogues in ergodic theory and Borel dynamics.

\begin{corollary} \label{approximation} Let $T$ be an aperiodic
homeomorphism and let $U(T) = U(T;\mu_1,...,\mu_n;\e)$ be a
$\tau$-neighborhood of $T$. Then there exists a periodic homeomorphism $P\in
\per_0$ such that $P\in U(T)\cap [T]$.
\end{corollary}
{\it Proof.} Apply Corollary \ref{C1} for $n=2$. We obtain a clopen
partition of $X$ into a finite number of towers $\Xi=\{\xi_1,\ldots,\xi_q\}$
such that for all $i=1,\ldots,k$
$$
\mu_i(\bigcup\limits_{j=1}^q\bigcup\limits_{l=0}^{h_j- 2}T^lB_j)>1-\e
$$
where $B_j$ and $h_j$ are the base and the height of $\xi_j$ respectively.
Now we can define the periodic homeomorphism $P$ as follows.
$$
Px=\left\{
\begin{array}{ll}
Tx, & \mbox{ if }x\in
\bigcup\limits_{j=1}^q\bigcup\limits_{l=0}^{h_j-2}T^lB_j,
\\
T^{-(h_j-1)}x & \mbox{ if }x\in T^{h_j-1}B_j\mbox{ for some }j= 1,...,q.
\end{array}\right.
$$
Obviously, $P\in [T]\cap \per_0$. It is easy to see that $P$ also belongs to
$U(T)$. \hfill$\square$

\begin{corollary}\label{nowhere} $\mathcal{A}p$ is a closed nowhere dense
subset in $(Homeo(X),\tau)$.
\end{corollary}
{\it Proof}. It was proved in [BDK2] that $\mathcal{A}p$ is closed in
$(\h,\tau)$. It follows from Corollary \ref{approximation} that
$\mathcal{A}p$ is a subset of $(Homeo(X),\tau)$ with empty interior.
\hfill$\square$

\begin{corollary}\label{RP} Let $S\in\ap$. Then, for any $R\in\ap$ and $U(R)
= U(R;\mu_1,\ldots,\mu_p;\e)$, there exists some $T\in Homeo(X)$
such that $T^{-1}ST\in U(R)$. In other words, $Ap$ is the
$\tau$-closure of $\{T^{-1}ST\,:\,T\in Homeo(X)\}$.
\end{corollary}
{\it Proof}. Given $R\in\ap$,\ $\mu_1,\ldots,\mu_p\in \M$, and $\e> 0$,
apply Corollary \ref{C1} with $n=3$. We get a clopen partition of $X$ into a
finite number of $R$-towers $\Xi=\{\xi_1,\ldots,\xi_q\}$ such that
$$
\mu_i(\bigcup\limits_{j=1}^q\bigcup\limits_{k=0}^{h_j-3}R^kB_j)>
1-\epsilon,\ \ i=1,2,\ldots,q,
$$
where $B_j$ is the base of $\xi_j$ and $h_j$ is its height.

Put $M=\max\{h_j\,:\,1\leq j\leq q\}$. Apply Proposition \ref{cover} with
$n=M$ to the homeomorphism $S\in \ap$ and find a clopen partition of $X$
into a finite number of $S$-towers $\Lambda= \{\lambda_1,\ldots,\lambda_l\}$
such that the height $n_i$ of each tower $\lambda_i$ is at least $M$. We can
refine either $\Xi$ or $\Lambda$ to get the same number of towers in $\Xi$
and $\Lambda$. Thus, without loss of generality, we can set $q=l$. Denote by
$Z_j$ the base of $\lambda_j$. Let $Q_j$ be a homeomorphism which maps $B_j$
onto $Z_j,\ j=1,...,q$, and let $Q$ be a homeomorphism which sends
$\bigcup\limits_{j=1}^qR^{h_j-1}B_j$ onto
$\bigcup\limits_{j=1}^q\bigcup\limits_{k=h_j-1}^{n_j-1} S^kZ_j$. Define the
homeomorphism $T$ as follows.
$$
Tx=\left\{
\begin{array}{ll}
S^iQ_jR^{-i}, & \mbox {if}\ x\in R^iB_j,\ 0\leq i\leq h_j-2,\ 1\leq j\leq
q\\
\\
Qx, & \mbox{ if }x\in \bigcup\limits_{j=1}^q R^{h_j-1}B_j
\end{array}\right.
$$

In such a way, the homeomorphism $T$ is defined on the set $X$. It is not
hard to see that $E_0(R,T^{-1}ST)=\{x \in X : Rx \neq
T^{-1}STx\}\subset\bigcup\limits_{j=1}^q\bigcup\limits_{i=1,2}R^{h_j-i}B_j$.
Therefore, $T^{-1}ST\in U(R)$. \hfill$\square$
\\

{\it Acknowledgement}. We would like to thank J.~Feldman,
A.~Kechris, J.~Kwiatkowski, and B.~Miller for helpful discussions.
The first named author thanks the University of New South Wales
for the warm hospitality and the Australian Research Council for
its support.
\\

\end{document}